\documentclass[12pt]{article}
\usepackage{amssymb}
\usepackage{amsmath}
\usepackage{latexsym}
\usepackage{amsthm}
\usepackage{color}

\theoremstyle{plain}
\newtheorem{theorem}{Theorem}%[section]
\newtheorem{lemma}{Lemma}%[section]
\theoremstyle{definition}
\begin{document}

\title{\textbf{\large Explicit estimates on the torus for the 
sup-norm and the crest factor of solutions of the Modified 
Kuramoto-Sivashinky Equation in one and two space dimensions       
}}  
\author{\textbf{Michele V. Bartuccelli$^\ast,$ Jonathan H. Deane$^\ast,$ 
Guido Gentile$^{\ast \ast}$}
\vspace{2mm}
\\ \small 
$^\ast$Department of Mathematics,
University of Surrey, Guildford, GU2 7XH, UK \\ 
\small $^{\ast \ast}$Dipartimento di Matematica e Fisica, 
Universit\`a Roma Tre, Roma, I-00146, Italy}

\maketitle

\begin{abstract}
\noindent 
We consider the Modified Kuramoto-Sivashinky Equation (MKSE) in one and two space dimensions
and we obtain explicit and accurate estimates of various Sobolev norms of the solutions.
In particular, by using the sharp constants which appear in the functional interpolation 
inequalities used in the analysis of partial differential equations, we evaluate explicitly the sup-norm of the solutions
of the MKSE. Furthermore we introduce and then compute the so-called 
crest factor associated with the above solutions. The crest factor provides information
on the distortion of the solution away from its space average and therefore, if it is large,
gives evidence of strong turbulence. Here we find that the time average of the crest factor
scales like $\lambda^{(2d-1)/8}$ for $\lambda$ large, where $\lambda$ is  the bifurcation parameter
of the source term and $d=1,2$ is the space dimension. This shows that strong turbulence
cannot be attained unless the bifurcation parameter is large enough.
\end{abstract}

\vspace*{0.2in}
\noindent{\bf Short title:} On the Crest Factor for the Modified 
Kuramoto-Sivashinsky Equation.

\vspace*{.1in}
\noindent {\bf Keywords:} Dissipative Partial Differential Equations; 
Interpolation Inequalities; 
Best Constants; Analysis of Solutions; Crest Factor.  

\vspace*{.1in}
\noindent{\bf Classification numbers:} 35B40, 35B45, 35G20, 35K25, 
46E20, 46E35

\newpage

\normalsize

%Creating a footnote is easy.\footnote{An example footnote.}
%\footnote{AMS Subject Classification: 46E20, 46E35; 35B40, 35B45.}

%%%%%%%%%%%%%%%%%%%%%%%%%%%%%%%%%%%%%%%%%%%%%%%%%%%%%%%%%%%%%%%%%%%%%%%%%
%%%%%%%%%%%%%%%%%%%%%%%%%%%%%%%%%%%%%%%%%%%%%%%%%%%%%%%%%%%%%%%%%%%%%%%%%
\section{\large Introduction}\label{sec1}
%%%%%%%%%%%%%%%%%%%%%%%%%%%%%%%%%%%%%%%%%%%%%%%%%%%%%%%%%%%%%%%%%%%%%%%%%
%%%%%%%%%%%%%%%%%%%%%%%%%%%%%%%%%%%%%%%%%%%%%%%%%%%%%%%%%%%%%%%%%%%%%%%%%

Accurate analysis of solutions of partial differential equations (PDEs) 
is an essential part in our understanding of many of the features of 
physical and biological phenomena. There are various approaches that
strive to obtain detailed information on the behaviour of solutions of PDEs.  
In this work we use functional analysis methods and we employ the latest explicit 
and sharp estimates for the embedding constants appearing in the functional inequalities
widely used in the study of any PDE. More precisely we have computed, as 
accurately as we possibly can, the estimates for some of the classical Sobolev 
norms of solutions of a model very close to some classical PDEs to which it reduces
in particular cases. In the following we will refer to our model as the 
Modified Kuramoto-Sivashinky Equation (MKSE); in two space dimensions it reads
\begin{equation} \label{MKS} 
   u_{t} = - \Delta^2 u - \Delta u + \lambda u -u^{3} - u(u_x + u_y) ,
\end{equation}
where $\Delta$ is the Laplacian, $u=u(x,y,t)$
for $(x,y)\in \Omega=[0,L]^2$, with $L > 0$ and $t>0$, subject to the initial 
condition $u(x,y,0)=u_0(x,y)$ and periodic boundary conditions 
on the boundary of $\Omega$. The real constant $\lambda$ is called the \emph{bifurcation parameter};
since we are mainly interested in the behaviour of the system for large $\lambda$,
for simplicity we take $\lambda>0$.
In this work we will obtain accurate estimates of some Sobolev norms of the MKSE
such as the the $L^\infty$ norm of its solutions. Furthermore we have introduced 
an important concept in the analysis of the behaviour of solutions of 
dissipative PDEs, namely
the so-called ``crest factor'', which is defined as the ratio between the $L^\infty$
and the $L^2$ norm of solutions. It has therefore the dimension of the square root 
of the inverse of the ``volume'' 
of the torus in $d$ spatial dimensions, and hence it can be made dimensionless by 
multiplying it  by $L^\frac{d}{2}.$   
The crest factor contains important informations
on the ``distortions'' between the amplitude and the 
$L^2$ norm of the solution. It is in fact a standard measurement used
in turbulence experiments in fluid dynamics.  
Effectively what it says is that if it is of order one then the dynamics is
relatively ``mild'', in the sense that the solution does not have major excursions 
in space-time. However, when the maximum amplitude of the solution becomes much 
larger with respect to its spatial average, then the solution does have strong
deviations in space and time; these strong intermittent fluctuations away from the
averages are one of the hallmarks of hard turbulence.  This phenomenon 
is now well established in many physical contexts such as, for example, in 
fluid convection. Thus the main aim of this work is to estimate in an explicit 
and accurate manner both some classical Sobolev norms of the solutions of the MKSE 
and the associated crest factor of these solutions.

Going back to our model first note that in the one space dimensional case  the (\ref{MKS}) naturally reduces
to the classical Kuramoto-Sivashinsky equation in the presence of a source term 
and a saturation term, namely one has
 \begin{equation} \label{MKS1d}
   u_{t} = - u_{xxxx} - u_{xx} + \lambda u -u^{3} - u u_x ,
\end{equation}
for $x \in \Omega =[0,L]$, with $L>0$, and $\lambda>0$.

Note also that by neglecting the last term in (\ref{MKS1d}) it reduces to another 
classical dissipative PDE, namely the Swift-Hohenberg equation. Both the 
Kuramoto-Sivashinsky equation and the Swift-Hohenberg equation have been 
extensively investigated because of their
fundamental importance in many mathematical, physical, biological and other 
contexts. So the literature on these two PDEs is huge and here we are forced 
to cite only a few of the relevant papers devoted to them: see for
example ~\cite{SH,CH,KT,S,KZ,HN,Tadmor,CFNT,JKT,CEES,PZ}.
  
The layout and main results of the paper are as follows:
in Section \ref{sec2} we state some standard functional setting and the notation
used in this work. In Section \ref{sec3} we obtain explicit and accurate estimates 
for the sup-norm of the solutions of the MKSE in one and two spatial 
dimensions. These estimates are stated after proving the 
Lemmas \ref{lemma1},\ref{lemma2},\ref{lemma3} 
and Theorem $1.$ In Section \ref{sec4} we compute the time averaged dissipative
length scale also in one and two spatial dimensions. Finally in Section
\ref{sec5} we obtain the ``crest factor'' of the solutions of the MKSE and we 
express the conclusion and open problems.  
  
\section{Functional Settings and Notation}  \label{sec2}  

Let us first give a brief standard preliminary functional setting and notation
~\cite{Adams,Robinson,Mazja,Temam}.  Denote by $\Omega = [0,L]^d$ the
$d-$dimensional
torus\,; for any scalar function $\phi(x)$ with $x \in \Omega$
let $\|\phi\|_p^p = \int_{\Omega} |\phi(x)|^p \,{\rm d}x$ be the norm associated with
the Banach space of
$\Omega-$periodic functions\,; we also define the $L^{\infty}$ norm as
\begin{equation} \nonumber %\label{norm1}
\|\phi\|_{\infty} = \sup_{x \in \Omega} |\phi(x)|\,.
\end{equation}
For $p=2$ we denote by $L^2(\Omega)$
the Hilbert space of $\Omega-$periodic functions $\phi$ with $\|\phi\|_{2}<+\infty$.
Given a multi-index $\vec{n}=(n_1,n_2,\ldots,n_d)$, with all the $n_i$ non-negative integers, let
$|\vec{n}|=n_1+\ldots+n_d$ and
\begin{equation} \nonumber %\label{Dn}
   D^{\vec{n}} := \frac{\partial^{|\vec{n}|}}
   {\partial x_{1}^{n_1} \partial x_{2}^{n_2} \cdot \cdot \cdot \partial x_{d}^{n_d}} ,
\end{equation}
and let
\begin{equation} \nonumber %\label{Hnze}
     {H}^{n} := \Bigl\{ \phi :
     \int_{\Omega} (D^{\vec{n}} \phi)^2 {\rm d}x < + \infty
     \hbox{ for all } \vec{n} \hbox{ such that } |\vec{n}|=n \Bigr\}   ,
\end{equation}
together with
\begin{equation}\label{Hnzero}
     \|\phi \|_{ {H}^{n}}^2 := \sum_{\substack{n_1,\ldots,n_d \ge 0 \\ n_1+ \ldots + n_d =n}}
      \frac{n!}{n_1! \cdot \cdot n_d!} 
      \|D^{\vec{n}} \phi\|_{2}^2 ,
\end{equation}
be the Sobolev space of $\Omega-$periodic functions with up to $n-$derivatives in $L^2(\Omega)$.
We also set $Du:=(\partial u /\partial x_1, \ldots, \partial u/\partial x_d)$. In (\ref{Hnzero}),
we naturally identify the functions having the same ``mixed'' partial derivatives, because it
is well known that the solutions of the MKSE are sufficiently smooth ~\cite{BV,Temam,Robinson}; 
for example we identify the differential operators 
\begin{equation}\label{Dniden}
   \frac{\partial^{n_1+n_2 + \cdot \cdot \cdot + n_d}}
   {\partial x_1^{n_1} \ldots \partial x_{i}^{n_i} \ldots \partial x_{j}^{n_j} \ldots \partial x_{d}^{n_d}} \,\equiv 
   \frac{\partial^{n_1+n_2 + \cdot \cdot \cdot + n_d}}
   {\partial x_1^{n_1} \ldots \partial x_{j}^{n_j} \ldots \partial x_{i}^{n_i} \ldots \partial x_{d}^{n_d}} ,
\end{equation}
and of course any other possible combination of the indices.  Also from Parseval's identity we have that
\begin{equation}\label{Parseval}
   \sum_{\substack{n_1,\ldots,n_d \ge 0 \\ n_1+ \ldots + n_d =n}}
      \frac{n!}{n_1! \ldots n_d!} 
      \|D^{\vec{n}} \phi\|_{2}^2 = 
    L^d \, \left(\frac{2 \pi}{L}\right)^{2n}
%    \sum_{\substack{n_1,\ldots,n_d \ge 0 \\ n_1+ \ldots + n_d =n}}
   \sum_{\vec{k} \in  \mathbb{Z}^d } |\vec{k}|^{2n}|\phi_{\vec{k}}|^2\,.  
\end{equation}
In (\ref{Parseval}) the Fourier series expansion has been used,
\begin{equation} \nonumber %\label{ex1}
\phi=\sum_{\vec{k}  \in \mathbb{Z}^d}\,
      \phi_{\vec{k}} \, e^{2 \pi i \vec{k} \cdot \vec{x} /L}\,,
\end{equation}
and
\begin{equation} \nonumber %\label{ex2}
   |\vec{k}|^2 = \vec{k} \cdot \vec{k}
    = k_1^2 + k_2^2 + \ldots + k_d^2 .
\end{equation}
By the same token the definition of Sobolev space can be extended to any {\it real} number $s$ as
\begin{equation}\label{Hreals}
    {H}^{s} \, = \Bigl\{  \phi=\sum_{\vec{k} \in \mathbb{Z}^d}
      \phi_{\vec{k}} \, e^{2 \pi i \vec{k} \cdot \vec{x} /L} :
      \overline{\phi}_{\vec{k}}=\phi_{-\vec{k}} \hbox{ and }
      \sum_{\vec{k} \in  \mathbb{Z}^d} |\vec{k}|^{2s}|\phi_{\vec{k}}|^2 
      < + \infty \Bigr\} ,
\end{equation}
and the corresponding norm is given by
\begin{equation} \nonumber
\| \phi \|_{H^s}^2 :=  L^d \, \left(\frac{2 \pi}{L}\right)^{2s}
%    \sum_{\substack{n_1,\ldots,n_d \ge 0 \\ n_1+ \ldots + n_d =n}}
   \sum_{\vec{k} \in  \mathbb{Z}^d } |\vec{k}|^{2s}|\phi_{\vec{k}}|^2\,.  
\end{equation}
These Sobolev spaces, defined on the $d-$dimensional torus,
are used below as we need to deal with the negative Laplacian
$A := - \Delta$ (as a self-adjoint unbounded operator) and its fractional powers.
More precisely, the eigenvalues of $A$
are given by the numbers $(2 \pi/L)^2 |\vec{k}|^{2}$, so the
domain of its powers $A^s$ is the set of functions such that
\begin{equation}\label{fractionalpowers}
   L^d \left(\frac{2 \pi}{L}\right)^{4s} \sum_{\vec{k} \in  \mathbb{Z}^d } |\vec{k}|^{4s}|\phi_{\vec{k}}|^2 
      = \| A^s \phi \|_2^2 < + \infty\,.   
\end{equation}
Thus in this paper, for any $s > 0$, we make the {\it formal} identification
\begin{equation} \nonumber %\label{identification}
   \|A^{\frac{s}{2}} \phi \|_2^2    =           
   \| (- \Delta)^{\frac{s}{2}} \phi \|_2^2 
     = L^d \left(\frac{2 \pi}{L}\right)^{2s} \sum_{\vec{k} \in  \mathbb{Z}^d} |\vec{k}|^{2s}|\phi_{\vec{k}}|^2,  
\end{equation}
\emph{provided it is understood that these operators are being used as
differential operators ``acting'' on functions in ${H}^{s}$,
according to (\ref{Hreals}) and (\ref{fractionalpowers}).}

%%%%%%%%%%%%%%%%%%%%%%%%%%%%%%%%%%%%%%%%%%%%%%%%%%%%%%%%%%%%%%%%%%%%%%%%
%%%%%%%%%%%%%%%%%%%%%%%%%%%%%%%%%%%%%%%%%%%%%%%%%%%%%%%%%%%%%%%%%%%%%%%%%

\section{Explicit Estimates of Sobolev norms of the MKSE}   \label{sec3}

\noindent In this section we wish to obtain explicit (and as accurately   
as we possibly can) estimates for various norm of solutions
of the MKSE. We then use such estimates to compute the corresponding 
crest factor associated to these solutions. In the light of this we then define
\begin{equation}\label{defH_n} 
    J_n \, := \,  \|u \|_{{H}^{n}}^2 
              = \sum_{\substack{n_1,\ldots,n_d \ge 0 \\ n_1+ \ldots + n_d =n}}
                \frac{n!}{n_1! \ldots n_d!} 
                \|D^{\vec{n}} u \|_{2}^2 .
\end{equation}
The MKSE has been defined in Section \ref{sec2} and it is given
by (\ref{MKS1d}) in $d=1$ and by (\ref{MKS}) in $d=2$,
in the domain $\Omega=[0,L]^d$, $d=1,2$, with $d$ being the spatial dimension.
The MKSE is known to have a unique solution for every initial
datum $u_0\in L^2(\Omega);$ the solution $u\in C([0,T];H)$, where  $H=L^2(\Omega)$,
for any $T>0;$ in addition the corresponding semigroup
$S_tu_0=u(t)$ has  a global attractor ${\cal A}\subset\!\subset H$
(for details see \cite{BV,Temam,Robinson}). Therefore all the calculations 
and estimates obtained below are {\it not formal,} but they 
reflect the actual behaviour of the solutions of the MKSE. Hence in the 
following we wish to find as accurately as possible estimates for the $J_n$ 
and then use them to obtain the corresponding estimates for the
$L^\infty$ norm of the solutions by using the sharp estimate found
in ~\cite{BG,BDZ,B} (see also ~\cite{X1,X2,Go,MS}.  

First note that one can show that the time-dependent functionals $J_n$ 
introduced above satisfy
a so-called {\it ladder differential inequality} \cite{BDGM,DG,BGO},
namely for any $n > d/2,$ where $d$ is the spatial
dimension, we have that
\begin{equation} \nonumber %\label{ladder}
      \frac{1}{2} \dot{J}_n \, \leq\, - J_{n+2} + J_{n+1} 
      + \lambda J_n  + \left( c_n \| u\|_{\infty}^2 + \tilde c_{n} \| Du \|_{\infty} \right) J_n,  
\end{equation}
where the constants $c_n$ and $\tilde c_n$ do not depend upon the solution $u=u(x,t)$.
Because we need to know explicitly all the constants appearing
in our analysis, we are somehow forced to restrict ourselves to the lower
values of the non-negative integer $n.$ In particular in the one-dimensional 
case we can restrict ourselves to the analysis of
$J_0$ and $J_1,$ which in $d=1$ are sufficient for having an upper bound
on the $\| u \|_{\infty}$ norm of the solution of any PDE. On the other
hand for the $d=2,3$ case we will have to analyze $J_2$ also.

Before starting our formal analysis let us make clear what we mean by the time-asymptotic
behavior of a given function of time $F(t)$. 
{\it From now on with an overbar over a given function of time $F(t),$ 
namely $\overline{F(t)},$ we mean the limit superior, taken
over all the initial conditions, as time goes to plus infinity. More
formally we mean that we are using the classical Gronwall inequality,
hence we take the limit superior as time goes to infinity 
and thence we consider the supremum over all the initial conditions.
Occasionally the set of initial conditions may be restricted
to the global attractor of the PDE under investigation, but this will be
clear from the context if not explicitly stated ~\cite{BV,Temam,Robinson}.} 

%\noindent {\bf
\subsection{Analysis in the one spatial-dimension case} 
%\\    
%\noindent

We can now start our analysis of our PDE on the torus in
one spatial dimension, namely we study
\begin{equation} \label{KSd1}
   u_{t} = - u_{xxxx} - u_{xx} + \lambda u -u^{3} - u u_x , 
\end{equation}
with periodic boundary conditions on $\Omega=[0,L]$.

In space dimension one it is sufficient to have control on the $J_0$
and the $J_1$ in order to have control on the sup norm of any solution of any
PDE. Thus we start with the analysis of $J_0(t)$.
%To simplify the notation, henceforth all the integrals 
%are meant as evaluated on the domain $\Omega=[0 , L]^d$, unless otherwise stated.

\begin{lemma}  \label{lemma1}
The time-asymptotic behaviour of $J_0(t),$ namely $\overline{J}_0,$ is given by 
\begin{equation}\label{H0lemma}
 \overline{J}_0 = \limsup_{t \to \infty} J_0(t) \leq 
 L \Bigl( \lambda + \frac{1}{4} \Bigr).  
\end{equation}
\end{lemma}

\noindent \textbf{Proof\,:} 
By taking the time-dependent quantity 
$J_0(t) = \int_{\Omega} u^2(x,t) \, {\rm d}x$ %:= \int u^2$
and differentiating it with respect to time one finds 
\begin{equation}\label{J0}
      \frac{1}{2} \dot{J}_0 = - J_{2} + J_{1} + \lambda J_0 
     - \int_{\Omega} (u)^4\,{\rm d}x.
\end{equation}
Note that the contribution from the last term in (\ref{KSd1}) is zero on periodic boundary
conditions. Also note that, for non-trivial behaviour one can see that we must have 
a restriction on the values of the parameter $\lambda;$ in fact,
after splitting the $J_{1}$ term by using first a Cauchy-Schwarz inequality
and then a Young inequality, namely
\begin{equation} \nonumber
J_1 \leq (J_2)^{\frac{1}{2}} (J_0)^{\frac{1}{2}} =
(2J_2)^{\frac{1}{2}} \Bigl( \frac{J_0}{2} \Bigr)^{\frac{1}{2}} 
\leq J_2 + \frac{1}{4} J_0 ,
\end{equation}
and also noting that $\displaystyle{- \int_{\Omega} (u)^4\,{\rm d}x
\leq - \frac{J_0^2}{L}},$ it follows that (\ref{J0}) becomes 
\begin{equation}\label{J0attr}
      \frac{1}{2} \dot{J}_0 \leq  \Bigl( \lambda +  \frac{1}{4} \Bigr) J_0 
     - \frac{J_0^2}{L} .
\end{equation}
Hence one can see that if $\lambda \leq - 1/4$ 
the zero solution becomes a global attractor. Since we have taken $\lambda>0$ we
are excluding such a situation. Thus going back to our analysis of $J_0$ we have to study (\ref{J0attr}).
By standard analysis one can see that the fixed points of the
corresponding nonlinear ordinary differential equation are given by
$J_0=0, L(\lambda + \frac{1}{4})$  
with $0$ being unstable and $L(\lambda + \frac{1}{4})$ being stable.
Thus the long-time asymptotic behaviour of $J_0$ (denoted with
$\overline{J}_0$) satisfies (\ref{H0lemma}). In particular
it is independent of the initial condition $u(x,t=0) = u_0(x).$  
\hfill$\blacksquare$ \\

We now turn our analysis to the estimate of $J_1$.

\begin{lemma} \label{lemma2}  
The time-asymptotic behaviour of $J_1(t),$ namely $\overline{J}_1,$
is given by
\begin{equation}\label{H1lemma}
  \overline{J}_1 := \limsup_{t \to \infty} J_1(t) \leq  
  \sqrt{\frac{24 \lambda + 13}{11}} L \Bigl( \lambda + \frac{1}{4} \Bigr).  
\end{equation}
\end{lemma}

\noindent \textbf{Proof\,:} Here we take the time-dependent quantity 
$J_1(t) = \int_{\Omega} (u_x(x,t))^2 {\rm d}x$ %:= \int (u_x)^2$
and differentiating it with respect to time we find
\begin{equation} \nonumber \label{J1}
      \frac{1}{2} \dot{J}_1= - J_{3} + J_2 + \lambda J_1 
     - 3 \int_{\Omega} u^2 (u_x)^2 \, {\rm d}x - \int_{\Omega} (u_x)^3 \, {\rm d}x
     - \int_{\Omega} (u)(u_x)(u_{xx}) \, {\rm d}x .
\end{equation}
An integration by parts on the last term gives
\begin{equation} \nonumber %\label{JJ1}
      \frac{1}{2} \dot{J}_1= - J_{3} + J_2 + \lambda J_1 
     - 3 \int_{\Omega} u^2 (u_x)^2 \, \, {\rm d}x - \int_{\Omega} (u_x)^3 \, {\rm d}x
     + \int_{\Omega} (u_x)^3 \, {\rm d}x + \int_{\Omega} (u)(u_x) (u_{xx}) \, {\rm d}x ;
\end{equation}
hence two terms cancel out and then by performing first a Cauchy-Schwarz
inequality and then a judicious Young inequality so as to generate the terms
$3 \int u^2 (u_x)^2 + \frac{J_2}{12}$ one obtains
\begin{equation} \nonumber %\label{J11}
      \frac{1}{2} \dot{J}_1 \le  - J_{3} + J_2 + \lambda J_1 + \frac{1}{12} J_2 .  
\end{equation}
By using a Young inequality on the term $J_2$ and simplifying we arrive at
\begin{equation}\label{JJ11}
      \dot{J}_1 \, \leq \,- \frac{11}{12} J_3 + \left( 2 \lambda + \frac{13}{12} \right) J_1.  
\end{equation}

We now use the inequality \cite{BDGM,DG,BGO} 
\begin{equation} \label{beauty} 
J_p \leq J_{p+r}^{\frac{q}{r+q}} J_{p-q}^{\frac{r}{r+q}} , \qquad p \ge q , \quad r \ge 0 ,
\end{equation}
with $p=1$, $r=2$ and $q=1$ to obtains $- J_3 \leq - J_1^3/J_0^2$. Hence inserting
this into (\ref{JJ11}), so as to obtain
\begin{equation}\label{J11rr}
      \dot{J}_1 \leq - \frac{J_1^3}{J_0^2} + 
       \left(\frac{24 \lambda + 13}{11}\right) J_1.
\end{equation}
and performing a similar analysis to that
used in obtaining the estimate (\ref{H0lemma}), one finds
\begin{equation} \nonumber %\label{J111}
  \overline{J}_1 := \limsup_{t \to \infty} J_1(t) \leq  
  \sqrt{ \frac{24 \lambda + 13}{11}} \, \overline{J}_0 ,
  \end{equation}
which, together with (\ref{H0lemma}), yields the result.
\hfill$\blacksquare$  \\

By using the estimates above it is interesting to obtain the
corresponding estimate for the $\|u \|_{\infty}$ of the
solution in the $d=1$ case. Here we can apply the sharp results found in \cite{BG,BDZ,B}:
for any function $u\in H^{1+\epsilon}$ one has 
\begin{equation}  \label{supnorm}
   \|u \|_{\infty} \, \leq \, \left(\frac{\zeta(1+ \epsilon)}{\pi} 
    \right)^{\frac{1}{2}}  
   \|(- \Delta)^{\frac{1+\epsilon}{4}} u \|_2\, +
    L^{- \frac{1}{2}} J_0^{\frac{1}{2}},  
\end{equation}
where $\epsilon > 0$ and 
\begin{equation} \label{RZfunc}
  \zeta(1+ \epsilon)= \sum_{n \geq 1} \frac{1}{n^{1+\epsilon}}
\end{equation}
is the Riemann zeta function.  The last term in (\ref{supnorm}) takes into account the mean of $u$.
By taking the value $\epsilon =1$ we therefore obtain
\begin{equation} \label{rriecat}
    \|u \|_{\infty} \, \leq \, 
     \sqrt{\frac{\pi}{6}} \,\| D u \|_{2} +
     L^{- \frac{1}{2}} J_0^{\frac{1}{2}} = 
     \sqrt{\frac{\pi}{6}} \, J_1^{\frac{1}{2}} +
     L^{- \frac{1}{2}} J_0^{\frac{1}{2}};   
\end{equation}
thus by using (\ref{H0lemma}) and (\ref{H1lemma}) we obtain
\begin{equation}  \label{supsol}
  \overline{\|u\|}_{\infty} \, \leq \, 
  \left(\frac{L \pi}{24}(4 \lambda + 1) 
  \sqrt{\frac{24 \lambda + 13}{11}}\right)^{\frac{1}{2}}  
  + \frac{\sqrt{4 \lambda + 1}}{2}.  
\end{equation}
%

%\noindent {\bf
\subsection{Analysis in the two spatial-dimensions case}
%\\ 

We can now turn our attention to the two-dimensional case having
domain $[0, L]^2$; as it is well 
known in this case having control on the $J_1$ norm alone is not
sufficient, but it is necessary to have control on the $J_2$ norm as well.
Before actually computing the time-asymptotic behaviour of $J_2$ we note 
that the estimates for $\overline{J}_0$ and $\overline{J}_1$ in two
spatial dimension are different because of the nonlinear terms;
indeed all we have to do is estimating the nonlinear part as best as 
we can. We start with the estimate of $J_0.$ Here the only difference 
with respect to the $d=1$ case comes from the term 
$\displaystyle{- \int  (u)^4\,{\rm d}x\,{\rm d}y
\leq - J_0^2/L^2}$; it follows that the differential inequality
for $J_0(t)$ becomes
\begin{equation} \nonumber %\label{J02dattr}
      \frac{1}{2} \dot{J}_0 \leq  \Bigl( \lambda +  \frac{1}{4} \Bigr) J_0 
     - \frac{J_0^2}{L^2}. 
\end{equation}
Therefore one obtains for the time-asymptotic behaviour of $J_0(t)$ the estimate 
\begin{equation} \label{J0asymin2d}
 \overline{J}_0 := \limsup_{t \to \infty} J_0(t) \leq 
 L^2 \Bigl( \lambda + \frac{1}{4} \Bigr).  
\end{equation}

Similarly for the time-asymptotic behaviour of $J_1$ one finds that
\begin{equation} \nonumber %\label{J2d}
 \frac{1}{2} \dot{J}_1= - J_{3} + J_2 + \lambda J_1 
 - \sum_{|\vec{n}|=1} \,\int_{\Omega} \left[ (D^{\vec{n}} u)D^{\vec{n}}(u^3)\,{\rm d}x\,{\rm d}y
 + \,(D^{\vec{n}}u)[D^{\vec{n}}(uu_x + uu_y)]\,{\rm d}x\,{\rm d}y \right] .
\end{equation}
Hence by neglecting the negative definite term given by the first summation 
and by expanding all the derivatives present in the second summation one arrives at
  \begin{equation} \label{J22d}
  \frac{1}{2} \dot{J}_1\le  - J_{3} + J_2 + \lambda J_1 
  + \sqrt{\frac{24}{\pi}} J_1 J_2^{\frac{1}{2}}. 
  \end{equation}  
Hence a similar analysis to the one done for obtaining the time-asymptotic 
behaviour of $J_0$ gives the estimate
\begin{equation}\label{J111in2d}
  \overline{J}_1 := \limsup_{t \to \infty} J_1(t) \leq  
  \left[\frac{5}{3} + \frac{5}{6} (4)^{\frac{14}{5}} 
  \Bigl(\frac{6}{\pi}\Bigr)^{\frac{3}{5}}\right]^{\frac{1}{3}}
  \overline{J}_0 \, \leq \, \left[\frac{5}{3} + \frac{5}{6} (4)^{\frac{14}{5}}
  \Bigl(\frac{6}{\pi}\Bigr)^{\frac{3}{5}}\right]^\frac{1}{3} L^2 \Bigl(\lambda + \frac{1}{4}\Bigr).  
\end{equation}
We now turn our attention to the analysis of $J_2(t);$ the
corresponding first order non-linear differential equation is given by
\begin{equation} \label{J2d}
 \frac{1}{2} \dot{J}_2= - J_{4} + J_3 + \lambda J_2 
 - \sum_{|\vec{n}|=2} \,\int_{\Omega} \Bigl( (D^{\vec{n}} u) D^{\vec{n}}(u^3)
 + \,(D^{\vec{n}}u) [D^{\vec{n}}(uu_x + uu_y)] \Bigr) {\rm d}x\,{\rm d}y .
\end{equation}
where the terms in the summations represent the non-linear terms. 
Their accurate estimates is given by the following  result.

\begin{lemma} \label{lemma3}
The nonlinear terms above obeys the estimate  
\begin{equation} \nonumber %\label{nonlin}
 - \sum_{|\vec{n}|=2} \,\int_{\Omega} \Bigl( (D^{\vec{n}} u)D^{\vec{n}}(u^3) 
 + (D^{\vec{n}}u)[D^{\vec{n}}(uu_x + uu_y)] \Bigr) {\rm d}x\,{\rm d}y 
\leq \, \frac{78}{\pi} J_{1} J_{2} +5 \|Du\|_\infty J_2.   
\end{equation}
\end{lemma}  

\noindent \textbf{Proof\,:} We first analyse the terms
\begin{equation} \nonumber
 - \sum_{|\vec{n}|=2} \,\int_{\Omega} (D^{\vec{n}} u)D^{\vec{n}}(u^3)\,{\rm d}x\,{\rm d}y .
\end{equation}
One starts by making the explicit differentiation, thereby obtaining
\begin{eqnarray} %\label{NLT}
     & & - \sum_{|\vec{n}|=2} \,\int_{\Omega} (D^{\vec{n}} u)D^{\vec{n}}(u^3)\,{\rm d}x\,{\rm d}y =
     -6 \int_{\Omega} u(u_x)^2 u_{xx} \nonumber \\
     & & - 3 \int_{\Omega} u^2 (u_{xx})^2\,{\rm d}x\,{\rm d}y  
     - 6 \int_{\Omega} u(u_y)^2 u_{yy}\,{\rm d}x\,{\rm d}y  
     - 3 \int_{\Omega} u^2 (u_{yy})^2 \,{\rm d}x\,{\rm d}y \nonumber \\
     & & - 6 \int_{\Omega} u^2(u_{xy})^2\,{\rm d}x\,{\rm d}y 
     - 12 \int_{\Omega} u u_x u_y u_{xy}\,{\rm d}x\,{\rm d}y; \nonumber
\end{eqnarray}
integrating by parts the first, the third and the last terms and then
rearranging we obtain
\begin{eqnarray} %\label{NLT1}
   & & - \sum_{|\vec{n}|=2} \,\int_{\Omega} (D^{\vec{n}} u)D^{\vec{n}}(u^3)\,{\rm d}x\,{\rm d}y =
     2 \int_{\Omega} (u_x)^4\,{\rm d}x\,{\rm d}y \nonumber \\& &
    - 3 \int_{\Omega} u^2 (u_{xx})^2\,{\rm d}x\,{\rm d}y 
     + 2 \int_{\Omega} (u_y)^4 \,{\rm d}x\,{\rm d}y 
     - 3 \int_{\Omega} u^2 (u_{yy})^2 \,{\rm d}x\,{\rm d}y \nonumber \\ && 
     -6 \int_{\Omega} u^2(u_{xy})^2\,{\rm d}x\,{\rm d}y
     + 6 \int_{\Omega} (u_x)^2(u_y)^2\,{\rm d}x\,{\rm d}y
     + 6 \int_{\Omega} u u_{xx}(u_y)^2\,{\rm d}x\,{\rm d}y . \nonumber
\end{eqnarray}
By splitting the last two terms by applying first a Cauchy-Schwarz
inequality and then a Young inequality we get
\begin{eqnarray} % \label{NLT11}
 & & - \sum_{|\vec{n}|=2} \,\int_{\Omega} (D^{\vec{n}} u)D^{\vec{n}}(u^3)\,{\rm d}x\,{\rm d}y =
   2 \int_{\Omega} (u_x)^4\,{\rm d}x\,{\rm d}y \nonumber \\ 
   & & 
    - 3 \int_{\Omega} u^2 (u_{xx})^2\,{\rm d}x\,{\rm d}y
     + 2 \int_{\Omega} (u_y)^4 \,{\rm d}x\,{\rm d}y - 3 \int_{\Omega} u^2 (u_{yy})^2 \,{\rm d}x\,{\rm d}y
     \nonumber \\ 
     && - 6 \int_{\Omega} u^2(u_{xy})^2\,{\rm d}x\,{\rm d}y  
     + 3 \int_{\Omega} (u_x)^4\,{\rm d}x\,{\rm d}y  + 3 \int_{\Omega} (u_y)^4\,{\rm d}x\,{\rm d}y \nonumber \\ && 
     + 3 \int_{\Omega} u^2(u_{xx})^2\,{\rm d}x\,{\rm d}y + 3 \int_{\Omega} (u_{y})^4\,{\rm d}x\,{\rm d}y .  \nonumber
\end{eqnarray}
Simplifying we finally obtain that the nonlinear term can be estimated as
\begin{equation} \label{NLT111}
     - \sum_{|\vec{n}|=2} \,\int_{\Omega} (D^{\vec{n}} u)D^{\vec{n}}(u^3)\,{\rm d}x\,{\rm d}y
    \leq \, 5 \int_{\Omega} (u_x)^4\,{\rm d}x\,{\rm d}y  + 8 \int_{\Omega} (u_{y})^4\,{\rm d}x\,{\rm d}y.
\end{equation}

In the two-dimensional case we can use an
improved version of the Ladyzhenskaya inequality \cite{Alex2}, namely
for any mean zero function $\phi(x,y)$ on the $2d$ torus we have 
the inequality
   $$\int_{\Omega} (\phi(x,y))^4 \,{\rm d}x\,{\rm d}y \leq 
     \frac{6}{\pi} \int_{\Omega} (\phi(x,y))^2 \,{\rm d}x\,{\rm d}y 
     \int_{\Omega} |\nabla \phi|^2 \,{\rm d}x\,{\rm d}y . $$
Hence we can estimate the term $5 \int_{\Omega} (u_x)^4\,{\rm d}x\,{\rm d}y$ in (\ref{NLT111}) as
$$5 \int_{\Omega} (u_x)^4\,{\rm d}x\,{\rm d}y \, \leq \, \frac{30}{\pi} 
\left(\int_{\Omega} (u_x)^2\,{\rm d}x\,{\rm d}y \,\right) \left(\int_{\Omega} (u_{xx}^2 
 + u_{xy}^2)\,{\rm d}x\,{\rm d}y \right)$$
 and similarly
$$8 \int_{\Omega} (u_y)^4\,{\rm d}x\,{\rm d}y \, \leq \, \frac{48}{\pi} 
\left(\int_{\Omega} (u_y)^2\,{\rm d}x\,{\rm d}y \,\right) \left(\int_{\Omega} (u_{yy}^2 
 + u_{xy}^2)\,{\rm d}x\,{\rm d}y \right).$$

By noting that $\int_{\Omega} (u_x)^2\,{\rm d}x\,{\rm d}y \, \leq \, J_{1}$, 
$\int_{\Omega} (u_y)^2\,{\rm d}x\,{\rm d}y \, \leq \, J_{1}$,
$\int_{\Omega} (u_{xx}^2 + u_{xy}^2)\,{\rm d}x\,{\rm d}y \, \leq J_{2}$ and 
$\int_{\Omega} (u_{yy}^2 + u_{xy}^2)\,{\rm d}x\,{\rm d}y \, \leq J_{2}$, 
we therefore obtain
\begin{equation}\label{NLTfinal}
     - \sum_{|\vec{n}|=2} \,\int_{\Omega} (D^{\vec{n}} u)D^{\vec{n}}(u^3)\,{\rm d}x\,{\rm d}y
         \leq \, \frac{78}{\pi} J_1 J_2.
\end{equation} 

We now turn to the other remaining nonlinear terms;
again we start be expressing them explicitly, namely
\begin{eqnarray} %\label{NLTT}
 & & -   \sum_{|\vec{n}|=2} \,\int_{\Omega} (D^{\vec{n}} u)D^{\vec{n}}(u u_{x} + u u_{y})\,{\rm d}x\,{\rm d}y =  \nonumber \\ & &
    - \int_\Omega u_{xx} [u u_x + u u_y]_{xx} 
 - \int_\Omega 2 u_{xy} [u u_x + u u_y]_{xy}  
 - \int_\Omega u_{yy} [u u_x + u u_y]_{yy} = \, \nonumber \\ & & 
 - \frac{5}{2} \int_{\Omega} u_x (u_{xx})^2\,{\rm d}x\,{\rm d}y  
     - \frac{1}{2}  \int_{\Omega} u_y (u_{xx})^2 \,{\rm d}x\,{\rm d}y  
     - 2 \int_{\Omega} u_x u_{xx} u_{xy} \,{\rm d}x\,{\rm d}y  \nonumber \\ & & 
      - 3 \int_{\Omega} u_x (u_{xy})^2\,{\rm d}x\,{\rm d}y  
     - 2 \int_{\Omega} u_y u_{xx} u_{xy} \,{\rm d}x\,{\rm d}y  
     - 3 \int_{\Omega} u_y (u_{xy})^{2} \,{\rm d}x\,{\rm d}y  \nonumber \\ & & 
 - 2 \int_{\Omega} u_x u_{xy} u_{yy} \,{\rm d}x\,{\rm d}y  
     - \frac{1}{2}  \int_{\Omega} u_x (u_{yy})^2 \,{\rm d}x\,{\rm d}y  
     - 2 \int_{\Omega} u_y u_{xy} u_{yy} \,{\rm d}x\,{\rm d}y  \nonumber \\ & & 
 - \frac{5}{2} \int_{\Omega} u_y (u_{yy})^2\,{\rm d}x\,{\rm d}y   ,\nonumber
\end{eqnarray}
where any term with three derivatives has first been integrated by parts
to move one derivative away to the remaining terms in the integral. 
All integrals are of the form 
$$ \int_{\Omega} u_x(u_{xx})^2 \,{\rm d}x\,{\rm d}y , \quad
\int_{\Omega} u_x u_{xx} u_{xy} \,{\rm d}x\,{\rm d}y , \quad
\int_{\Omega} u_x (u_{xy})^2 \,{\rm d}x\,{\rm d}y ,$$
$$ \int_{\Omega} u_x u_{xy} u_{yy} \,{\rm d}x\,{\rm d}y , \quad
\int_{\Omega} u_x (u_{yy})^2 \,{\rm d}x\,{\rm d}y , $$
or with the variables $x$ and $y$ exchanged. We pull the terms $u_x$ or $u_y$
in the $L^\infty$ norm thereby obtaining, for instance,
$\int u_x(u_{xx}^2) \leq \|u_x\|_\infty J_{2,x}$, where with $J_{2,x}$ we mean the ``component  of $J_2$ along
the $x$ coordinate"; the other similar terms such as $\int u_y (u_{xx}^2)$,
$\int u_x(u_{yy}^2)$, etc. are handled in the same way. 
Other terms of the form, say, $\int u_x u_{xy} u_{yy}$ are dealt with by first pulling
out the $u_x$ term in $L^\infty,$ then applying a 
Cauchy-Scharwz to the two remaining terms
and then splitting the two terms with a Young inequality.
We collect all the terms together thereby finally obtaining
\begin{equation} \nonumber %\label{NLTTT}
-  \sum_{|\vec{n}|=2} \,\int_{\Omega} (D^{\vec{n}} u)D^{\vec{n}}(u u_{x} + u u_{y})\,{\rm d}x\,{\rm d}y 
   \leq \, \|Du\|_\infty (5 J_{2,xx} + 2 \cdot 5 J_{2,xy} + 5 J_{2,yy})
   \leq 5 \|Du \|_\infty J_2,
\end{equation} 
where we have used that $\|u_x\|_{\infty},\|u_y\|_{\infty} \le \|Du\|_\infty$.
The last estimate, together with (\ref{NLTfinal}) implies the result.
\hfill$\blacksquare$ \\

By using the results obtained above we can now prove the following result.

\begin{theorem} \label{theorem1}
The time-asymptotic behaviour of $J_2(t),$ namely $\overline{J}_2,$ satisfies
\begin{eqnarray} %\label{H2asymin2d}
  \overline{J}_2 & \!\!\!\!  \leq \!\!\!\! &
   \overline{J}_{0}^{\frac{3}{2}} 
    \left[108 + 4 \lambda^2 + 108 \left(\frac{5}{\sqrt{\pi}}\right)^4 \overline{J}_{0}^2 +
    108 \left(\frac{78}{\pi}\right)^4 \overline{J}_{0}^4 \right]^\frac{1}{2} \nonumber \\
  & \!\!\!\!  \leq \!\!\!\! & 
    \left[ L^3 \Bigl(\frac{4 \lambda + 1}{4}\Bigr)^{3} \Bigl(
    108 + 4 \lambda^2 + 108 L^2 \left(\frac{5}{\sqrt{\pi}}\right)^4 \Bigl(\frac{4 \lambda + 1}{4}\Bigr)^2 +
    108 L^4 \left(\frac{78}{\pi}\right)^4 \Bigl(\frac{4 \lambda + 1}{4}\Bigr)^4 \right]^\frac{1}{2} \!\!\! . \nonumber
\end{eqnarray}
\end{theorem}  

\noindent \textbf{Proof\,:} First we write the estimate for the time
derivative of $J_2$, as obtained from (\ref{J2d}) and Lemma \ref{lemma3}, namely 
\begin{equation} \label{J22}
      \frac{1}{2} \dot{J}_2 \, \leq \, - J_{4} + J_3 + \lambda J_2 
     + \frac{78}{\pi} J_{1} J_{2} + 5 \|Du\|_\infty J_2.  
\end{equation}
To handle the last term we use the (almost sharp) estimate 
$\| Du \|_\infty \leq \, \frac{1}{\sqrt{\pi}} J_{3}^{\frac{1}{4}} J_{1}^{\frac{1}{4}}$ \cite{IT}.
To absorb the (\ref{J22}) term we split it as follows:
\begin{eqnarray}
& & J_3 \leq J_{4}^{\frac{3}{4}} J_{0}^{\frac{1}{4}}   \leq \, \frac{1}{8} J_4 + 54 J_0 , \nonumber \\
& & \lambda J_2 \leq \lambda J_4^{\frac{1}{4}} J_0^{\frac{1}{2}} \leq \frac{J_4}{8} + 2 \lambda^2 J_0 , \nonumber \\
& & \frac{78}{\pi} J_1 J_2 \leq 
\frac{78}{\pi} \left( J_4^{\frac{1}{4}} J_0^{\frac{3}{4}}  \right) \left( J_4^{\frac{1}{2}} J_0^{\frac{1}{2}} \right) \leq 
\frac{78}{\pi} J_4^{\frac{3}{4}} J_0^{\frac{5}{4}}  \leq 
\frac{J_4}{8} + \frac{1}{4} \left( 216 J_{0}^5 \Bigl(\frac{78}{\pi} \Bigr)^4 \right) , \nonumber \\
& & 5 \|Du\|_\infty J_2 \leq  
\frac{5}{\sqrt{\pi}} J_3^{\frac{1}{4}} J_1^{\frac{1}{4}} J_2 \leq
\frac{5}{\sqrt{\pi}} J_4^{\frac{3}{4}} J_0^{\frac{3}{4}} \leq
\frac{J_4}{8} + \frac{1}{4} \left( 216 \Bigl( \frac{5}{\sqrt{\pi}} \Bigr)^4 J_{0}^3 \right) , \nonumber
\end{eqnarray}
where (\ref{beauty}) has bee used repeatedly. Using all of this one arrives at
  \begin{equation} \nonumber %\label{J22TRIS}
      \frac{\dot{J}_2}{2} \leq \, - \frac{J_4}{2} + \frac{1}{2} J_0 
      \left[108 + 4 \lambda^2 + 108 \left(\frac{5}{\sqrt{\pi}}\right)^4 J_{0}^2 + 108 \left(\frac{78}{\pi}\right)^4 J_{0}^4 \right] .     
  \end{equation}
Therefore the time-asymptotic behaviour of $J_2$ is given by 
  \begin{equation} \label{J222asymp}
  \overline{J}_2\, \leq \, \overline{J}_{0}^\frac{3}{2} 
   \left[108 + 4 \lambda^2 + 
      108 \left(\frac{5}{\sqrt{\pi}}\right)^4 \overline{J}_{0}^2 + 108 \left(\frac{78}{\pi}\right)^4 \overline{J}_{0}^4 \right]^\frac{1}{2} .  
  \end{equation}
By substituting the estimate for $\overline{J}_{0}$ we finally obtain the result.
\hfill$\blacksquare$  \\

Thus for the estimate of $\|u\|_{\infty}$ we use the result proved in
~\cite{BG}, where it is shown that 
on the two-dimensional torus $\Omega=[0,L]^2$, for every $\epsilon > 0$,
the $L^{\infty}$ norm of a mean zero scalar function 
$u \in {H}^{1+\epsilon}$ satisfies the estimate
\begin{equation}  \label{supres}
   \|\phi\|_{\infty} \, \leq \, \left[4 \zeta(1+ \epsilon) 
   \beta(1+ \epsilon) \right]^{\frac{1}{2}}
    L^{-1} \left(\frac{L}{2 \pi}\right)^{(1+\epsilon)} 
   \|(- \Delta)^{\frac{1+\epsilon}{2}} \phi\|_2\,,
\end{equation}
where the coefficient $4 \zeta(1+ \epsilon) \beta(1+ \epsilon)$ 
is sharp, and where
\begin{equation} \nonumber %\label{RBfns}
\zeta(1+ \epsilon)= \sum_{n \geq 1} \frac{1}{n^{1+\epsilon}}\,,
\qquad  
\beta(1+ \epsilon)= \sum_{n \geq 0} 
\frac{(-1)^n}{(2n + 1)^{1+\epsilon}}\,,
\end{equation}
are the Riemann zeta-function and Dirichlet series respectively.
Thus for the estimate of $\overline{\|u\|}_{\infty}$ we use
(\ref{supres})  with $\epsilon = 1,$ namely
\begin{equation} \nonumber \label{sssupsolution}
  \overline{\|u\|}_{\infty} \, \leq \, 
  \frac{L}{2 \pi^2} (\zeta(2) \beta(2))^{\frac{1}{2}} \, 
  \|\Delta u\|_2 \, + L^{-1} \overline{J}_0^{\frac{1}{2}} \leq
  \frac{L}{2 \pi^2} (\zeta(2) \beta(2))^{\frac{1}{2}} \, 
  \overline{J}_2^{\frac{1}{2}} + L^{-1} \overline{J}_0^{\frac{1}{2}}.
\end{equation}
By using the values for
$\zeta(2) \beta(2) = 6 \pi^{-2} K$ with $K=0.915965594...$ we obtain  
\begin{equation}  \nonumber %\label{ssssupsolution}
  \overline{\|u \|}_{\infty} \, \leq \,
   \frac{L}{2 \pi^3} \sqrt{6 K} \overline{J}_2^{\frac{1}{2}} +
   L^{-1} \overline{J}_0^{\frac{1}{2}},     
\end{equation}
where the estimate for $\overline{J}_2^{\frac{1}{2}}$ is provided by
(\ref{J222asymp}) and that for $\overline{J}_0^{\frac{1}{2}}$ is
provided by (\ref{J0asymin2d}).      

\section{The Crest Factor of Solutions of Dissipative PDEs} \label{sec4}   

So far we have obtained various Sobolev norms estimates of solutions 
of our equation, such us the estimates for $J_0, ~ J_1,~J_2$ and the 
corresponding estimate for the sup-norm. An important question which
naturally arise from our analysis is to investigate the so-called 
{\it crest factor} (also known as the {\it peak to average ratio}), namely 
the ratio between the $L^\infty$ norm and the $L^2$ norm of the solution:
\begin{equation} \label{crest}
    C_f \, := L^\frac{d}{2} \frac{\|u\|_{\infty}}{J_0^{\frac{1}{2}}}.   
\end{equation}
It is therefore by definition dimensionless and it contains important 
information on the ``distortions'' between the sup-norm (the amplitude)
and the $L^2$ norm of the solution. It is in fact a standard measurement 
used in turbulence experiments in fluid dynamics. 
The ideal result would be to have a time-pointwise estimate of $C_f.$ However
this is very difficult due essentially to the non-linearity of the equation.
Alternatively one could try to estimate the time-asymptotic behaviour of $C_f$, 
but this also proves to be very hard to handle and it is essentially due to the 
lack of knowledge of a ``decent'' lower bound on the quantity $J_0,$ namely an 
estimate of the form 
$J_0(t) \geq \alpha > 0.$ The problem of estimating the lower bound 
appears in many contexts
in the theory of nonlinear dissipative PDEs, such as for example in the theory
of the Navier-Stokes equations where it is notoriously very hard to find a
``proper'' lower bound for the energy even on the torus \cite{DFJ}.  
So in this work we will compute the time-average of the quotient between 
the $L^\infty$ norm and the $L^2$ norm of the solution, namely  
$\displaystyle{\big< \|u\|_{\infty}/J_0^{\frac{1}{2}} \big>.}$ 
First of all let
us derive sharp estimates for the $\|u\|_{\infty}$ of typical solutions $u(x,t).$
Note that in general we cannot assume that the solutions of our equation
have zero-mean. Hence we have to ``carry along'' the mean value of our solutions. Thus
define $u^{*}(t):=\int_\Omega u(x,t)\,{\rm d}x$ and write $u(x,t)=u^*(t)+u^{\prime}(x,t)$,
where $\int_\Omega u^{\prime}(x,t)\,{\rm d}x=0$. Then using the inequality
\begin{equation}\label{mean}
|u^*|=L^{-d}\left|\int_\Omega u(x)\,{\rm d}x\right| \leq 
L^{-\frac{d}{2}} J_0^{\frac{1}{2}} 
\end{equation}
and defining $J_{0}':=\|u'\|_{2}^2$, we obtain \cite{BGI}
\begin{equation}\label{A4}
\|u\|_\infty\le|u^*|+\|u^{\prime}\|_\infty\le
L^{-\frac{d}{2}} J_0^\frac{1}{2} +
c(n) (J_0')^\frac{2n-d}{4n} J_n^\frac{d}{4n}\,.
\end{equation}
with $n > 1/2$ and $c(n)$ a suitable constant,
where we have used a Gagliardo-Nirenberg inequality to obtain the estimate on $\|u^{\prime}\|_\infty$. 
By substituting $u=1$ in (\ref{mean}) we see that the
constant $L^{-\frac{d}{2}}$ is sharp.
Therefore we obtain the following estimate 
\begin{equation} \nonumber %\label{supwithav}   
   \frac{\|u\|_{\infty}}{J_0^{\frac{1}{2}}} \leq \frac{|u^*| + 
   \|u^{\prime}\|_\infty }{J_0^{\frac{1}{2}}} \leq L^{-{\frac{d}{2}}}   
   + \frac{\|u^{\prime}\|_\infty }{J_0^{\frac{1}{2}}}.
\end{equation}
Hence by using (\ref{A4}) we obtain
\begin{equation} \nonumber
\frac{ \| u \|_\infty }{J_0^{\frac{1}{2} } }\leq 
L^{-{\frac{d}{2}}} + c(n) \left(\frac{J_n}{J_0}\right)^\frac{d}{4n} 
\left(\frac{J_0'}{J_0}\right)^\frac{2n-d}{4n} .
\end{equation}
Thus our estimate for the crest factor is obtained by taking the time-average
\begin{equation} \nonumber %\label{crest1}
    \widetilde{C}_f := \left< L^{\frac{d}{2}} \frac{\|u\|_{\infty}}{J_0^{\frac{1}{2}}} \right>.   
\end{equation}
It is useful to concentrate on the ``pure'' distortion between the sup-norm 
and the $L^2$ norm for non-constant solutions (note that of course 
constant functions have crest factor equal to $1$). Bearing this in mind 
one obtains 
  \begin{equation} \label{crestav}
   \widetilde{C}_f = 1 + \overline{C}_{f} , \qquad \overline{C}_{f} := \left< L^{\frac{d}{2}} \frac{\|u^{\prime}\|_\infty }
   {J_0^{\frac{1}{2}}} \right> \leq c(n) L^{\frac{d}{2}} 
   \left< \left(\frac{J_n}{J_0}\right)^\frac{d}{4n} \right> ,
   \end{equation}  
where the last bound follows noting that $J_0' \le J_0$.
Note that, since one has trivially $\overline{C}_{f}=0$ if $u(x,t)$ does not depend on $x$,
in order to estimate the crest factor we may assume in the following that $u'\neq0$. Hence $J_n>0$ for all $n\ge 0$.

%\noindent {\bf
\subsection{Time-averaged crest factor in one spatial dimension}
% \\

In one spatial dimension it is sufficient to take $n=1$ in (\ref{crestav})  and so one has    
\begin{equation} \label{dll} 
\overline{C}_f  \le c(1)
\left< L^{\frac{1}{2}} \left(\frac{J_1}{J_0}\right)^\frac{1}{4} \right> .
\end{equation}
From \cite{Alex1} or Appendix A in \cite{BGI} we have that $c(1) =1.$
Thus one needs to derive as best as possible the time average 
of the quantity $\left(J_1/J_0\right)^\frac{1}{4}$.  
This is achieved as follows. First take the differential 
inequality (\ref{J11rr}) and divide throughout by $J_1$. This leads to
\begin{equation} \nonumber %\label{J11rradd}
 \frac{\dot{J}_1}{J_1} \leq - \left(\frac{J_1}{J_0}\right)^2 +
   \left(\frac{24 \lambda + 13}{11}\right).     
\end{equation}
Then we take the time average of both sides
of the inequality thereby getting   
\begin{equation} \nonumber %\label{ta}
      \left< \left(\frac{J_1}{J_0}\right)^2 \right> \, \leq \, 
      \left(\frac{24 \lambda + 13}{11}\right) ,
\end{equation}
where we have used that $J_1$ is bounded both from below and from above by two positive constants.
Going back to (\ref{dll}) one obtains (with $c^2(1) = 1$), 
\begin{equation} \label{dlll} 
  \overline{C}_f \le L^\frac{1}{2} \left< \left(\frac{J_1}{J_0}\right)^\frac{1}{4} \right>
  \leq L^\frac{1}{2} \left< \left(\frac{J_1}{J_0}\right)^2 \right>^\frac{1}{8}   
  \leq L^\frac{1}{2} \left(\frac{24 \lambda + 13}{11}\right)^\frac{1}{8} , 
  \end{equation}
which shows that $C_f = O(\lambda^{\frac{1}{8}})$ for large $\lambda$.

%\noindent {\bf
\subsection{Time-averaged crest factor in two spatial dimensions} 
%\\  

The strategy for obtaining the time-averaged crest factor in two 
spatial dimensions is similar to the one-dimensional case with the 
corresponding changes, namely here $d=2$ and also one has to 
insert the explicit values of the constants $c(2)$ 
in (\ref{crestav}). Also it is well known that in two spatial dimensions it is
sufficient to take $n=2,$ and so we need to estimate the quantity
\begin{equation} \label{Cf2}
\overline{C}_f \le  c(2) L \left< \left(\frac{J_2}{J_0}\right)^\frac{1}{4} \right>
\le  c(2) L \left< \frac{J_2}{J_0} \right>^{\frac{1}{4}} .
\end{equation}

So we start from the differential inequality (see (\ref{J22}))
\begin{equation} \nonumber %\label{J22BISr}
      \frac{1}{2} \dot{J}_2 \, \leq \, - \frac{1}{2} J_{4} + 
      \frac{1}{2} \left(2 \lambda + 1 + \frac{156}{\pi} J_1 
      +10 \|Du\|_\infty \right) J_2.
\end{equation}
We now use again the inequality
$\|Du\|_\infty \leq \, \frac{1}{\sqrt{\pi}} (J_3 J_1)^\frac{1}{4}$ and also (\ref{beauty}),
with $p=q=r=2$, and so we obtain
\begin{equation} \nonumber %\label{J22BISrr}
      \dot{J}_2 \, \leq \, - \frac{J_2^2}{J_0}
      + \left(2 \lambda + 1 + \frac{156}{\pi} J_1 + \frac{10}{\sqrt{\pi}} 
      (J_3 J_1)^\frac{1}{4} \right) J_2. 
\end{equation}
Similarly to the one-dimensional case we divide throughout by $J_2$ and 
then we take the time average of both sides of the inequality obtaining
\begin{equation} \nonumber %\label{taJ2/J0}
      \left<\frac{J_2}{J_0}\right> \leq (2 \lambda + 1) + 
       \frac{156}{\pi} \big< J_1 \big> + \frac{10}{\sqrt{\pi}} 
      \big< J_3 \big>^\frac{1}{4} \big< J_1 \big>^\frac{1}{4}, 
\end{equation}
where we have used the properties of the time average 
in order to obtain the last term. 
In order to estimate $\big< J_1 \big>$ we use (\ref{J0}), that we re-write here:
\begin{equation} \nonumber %\label{J0r}
      \frac{1}{2} \dot{J}_0 = - J_{2} + J_{1} + \lambda J_0 
     - \int_{\Omega} (u)^4\,{\rm d}x.
\end{equation}
By using $J_1 \leq J_2^{\frac{1}{2}} J_0^{\frac{1}{2}}$ and then splitting 
the right hand side with the Young inequality we obtain
\begin{equation} \nonumber %\label{J0rr}
      \frac{1}{2} \dot{J}_0 = - \frac{J_2}{2}  + 
      \Bigl(\lambda + \frac{1}{2}\Bigr) J_0,
\end{equation}
where here we have neglected the last term.
Time averaging both sides we finally get
\begin{equation} \label{J2av} 
   \big< J_2 \big> \leq (2 \lambda + 1) \big< J_0 \big> \leq \, 
   (2 \lambda + 1) \overline{J}_0 \leq \, 
   (2 \lambda + 1) L^2\Bigl(\lambda + \frac{1}{4}\Bigr).  
\end{equation}
Thus by first time averaging the inequality 
$J_1 \leq J_2^{\frac{1}{2}} J_0^{\frac{1}{2}}$ and then splitting the
time average of the product on the right hand side one obtains
\begin{equation}  \label{J1av} 
   \big< J_1 \big> \leq \big< J_2 \big>^\frac{1}{2} 
   \big< J_0 \big>^\frac{1}{2} \leq (2 \lambda + 1)^\frac{1}{2} 
   \big< J_0 \big>^\frac{1}{2} \big< J_0 \big>^\frac{1}{2}
   \leq (2 \lambda + 1)^\frac{1}{2} \big< J_0 \big>
   \leq (2 \lambda + 1)^\frac{1}{2} \overline{J}_0.    
\end{equation}
%In the last inequality we have used the result that in the global 
%attractor it holds that $\big< J_0 \big> \leq \overline{J}_0.$
We now estimate the other term, namely $\big< J_3 \big>.$
Here we use the time-average of the formula  (\ref{J22d}), obtaining
\begin{equation} \nonumber % \label{J3av} 
   \big< J_3 \big> \leq \big< J_2 \big> + \lambda  \big< J_1 \big> +
   \frac{24}{\pi} \big< J_1 J_2^\frac{1}{2} \big> 
   \leq \, \big< J_2 \big> + \lambda  \big< J_1 \big> +
   \frac{24}{\pi} \big< J_0^\frac{1}{2}  J_2 \big>
   \leq \, \big< J_2 \big> + \lambda  \big< J_1 \big> +
   \frac{24}{\pi} \big< J_2 \big> \overline{J}_0^\frac{1}{2}.
\end{equation}  
Therefore by inserting the estimates for $\overline{J}_0$,
$\big< J_1 \big>$ and $\big< J_2 \big>$ given by (\ref{J0asymin2d}),
(\ref{J1av}) and (\ref{J2av}), respectively, one finds
\begin{equation}  \label{J3av-estim} 
   \big< J_3 \big> \leq (2 \lambda + 1)^\frac{1}{2} L^2 \Bigl(\lambda + \frac{1}{4}\Bigr)
   \left[\lambda + (2 \lambda + 1)^\frac{1}{2} \left(1+ \sqrt{\frac{24}{\pi}}
   L\Bigl(\lambda + \frac{1}{4}\Bigr)^\frac{1}{2} \right) \right].
\end{equation}  
So going back to the computation of the crest factor in the space two dimensional case, in (\ref{Cf2}) 
we have to insert the value of the constant $c(2)$, which is $c(2) = \sqrt{1/\pi}$ \cite{IT},
and bound $\big< (J_2/J_0)^\frac{1}{4}\big>$ by using
the estimates (\ref{J1av}) and (\ref{J3av-estim}) found above for $\big< J_1 \big>$ and $\big< J_3 \big>$.
In particular, one finds $C_f = O(\lambda^{\frac{3}{8}})$ for $\lambda$ large.

\section{Conclusions and Open Problems} \label{sec5}  

In this work we have analysed various Sobolev norms of solutions of a
modified version of 
KSE, with the aim to estimate as accurately as possible both the sup-norm
of solutions and then the corresponding {\it crest factor.} 
More specifically, by using the {\it best} available 
explicit estimates for the coefficients which appear in the Sobolev norms
used, we have first derived explicit estimates for the $\overline{J}_0, 
\overline{J}_1, \overline{J}_2,$ namely their time-asymptotic behaviour,
and then we have used these estimates to compute the time-asymptotic 
behaviour of the $L^\infty$ norm of the solution, namely the
$\overline{\|u\|}_{\infty}$ in one and two space dimension.
We then addressed another very important indicator of the dynamics
of solutions of dissipative PDEs, namely the accurate estimate of the 
so-called \emph{crest factor}. This is defined as the ratio
between the $L^\infty$ norm of the solution and the $L^2$ norm of the solution:
\begin{equation} \label{crest}
C_f \, := L^{\frac{d}{2}} \frac{\|u\|_{\infty}}{J_0^{\frac{1}{2}}},   
\end{equation}
where $d$ is the spatial dimension. 
It is therefore a dimensionless pure number and it contains important 
information on the ``distortion'' between the ``amplitude`` and the 
$L^2$ norm of the solution. It is in fact a standard measurement used
in turbulence experiments in fluid dynamics.   

Let us now discuss the implications of the estimates we have found in both 
one and two space dimensions. In space dimension one we found that th time-average of $C_f$ is
\begin{equation} \label{crestt-1d} 
\widetilde{C}_f = 1 + O(\lambda^{\frac{1}{8}}) ,
\end{equation}
while in space dimension two we found
\begin{equation} \label{crestt-2d} 
\widetilde{C}_f =  1 + O(\lambda^{\frac{3}{8}}) .
\end{equation}

The two formulas above reveal some of the features related to the
dynamics of the solutions of our PDE. In fact in one space dimension
the time average of the ratio between the ``peak to the root mean square'' 
(the crest factor) scales
like (\ref{crestt-1d}) as a function of the positive parameter $\lambda.$
So for small $\lambda$ the distortion $\overline{C}_f=\widetilde{C}_f-1$ is small (as it should),
but what it really says  is that our PDE cannot have {\it major} excursions in space-time as $\lambda$ 
increases because the crest factor goes like $\lambda^\frac{1}{8},$
and $\frac{1}{8}$ is ``pretty small``. On the other hand as a function of the
parameter $L$ (the length of the torus) it scales like $\sqrt{L};$ this
shows that the crest factor are more sensitive to the length of the torus than to the parameter $\lambda$.

In the two space dimension case the crest factor shows (of course) 
stronger potential fluctuations. Indeed it scales like the 
$\lambda^\frac{3}{8}$ for large $\lambda$, which is naturally
much larger then in the one space dimension case. As a function of
$L$, for large $L$ it scales like $L^\frac{3}{2}$,which again naturally
is much larger than in the one space dimension case where it goes like $\sqrt{L}$ for large $L.$ 

It would be interesting to compute the crest factor for other important PDEs,
such as the Complex Ginzburg-Landau equation and the Navier-Stokes equations. As one can
infer from our analysis above, the crest factor sheds some light on the nature of the solutions
of any PDE. In particular, as a function of the parameters and the length of the torus, it gives 
important indication on the fluctuations of solutions away from their spatial average. Thus it
can provide insight on regimes of ``soft" and ``hard" turbulent behaviour of the solutions
of any dissipative PDE.

\vspace{.5cm}

\noindent{\bf Acknowledgements}:
\noindent It is a pleasure to acknowledge very helpful discussions with
Paolo Secchi and Davide Catania on ideas and techniques closely related to this work.


\begin{thebibliography}{99}

\bibitem{Adams} R A Adams,
{\it Sobolev spaces}, Academic Press, New York, 1975.
%
\bibitem{BV} A V Babin, M I Vishik,
{\it Attractors for evolution equations},
 Nauka, Moscow, 1989.
%
\bibitem{B} M V Bartuccelli,
{\it Sharp constants for the $L^{\infty}$-norm on the torus and applications to dissipative partial differential equations},
 Differential Integral Equations {\bf 27} (2014), no. 1-2, 59-80.
%
\bibitem{BB} M V Bartuccelli,
{\it Explicit estimates on the torus for the sup-norm and the dissipative length scale of solutions of the
Swift-Hohenberg equation in one and two space dimensions},
J. Math. Anal. Appl. {\bf 411} (2014), no. 1, 166-176.
%
\bibitem{BDGM} M V Bartuccelli, C R Doering, J D Gibbon, S A Malham,
{\it Length Scales in Solutions of the Navier-Stokes equations},
Nonlinearity {\bf 6} (1993), 549-568.
%
\bibitem{BDZ} M V Bartuccelli, J H B Deane, S Zelik,
{\it Asymptotic expansions and extremals for the critical Sobolev and Gagliardo-Nirenberg inequalities on a torus},
Proc. Roy. Soc. Edinburgh Sect. A {\bf 143} (2013), no. 3, 445-482.
%
\bibitem{BG} M V Bartuccelli, J D Gibbon,
{\it Sharp constants in the Sobolev embedding theorem and a derivation of the Brezis-Gallouet interpolation inequality},
J. Math. Phys. {\bf 52} (2011), no. 9, 093706, 9 pages.  
%
\bibitem{BGO} M V Bartuccelli, J D Gibbon, M Oliver,
{\it Length Scales in Solutions of the complex Ginzburg-Landau equation}, 
Phys. D {\bf 89} (1996), no. 3-4, 267-286.
%
\bibitem{BGI}  M V Bartuccelli, S A Gourley, A A Ilyin,
{\it Positivity and the attractor dimension in a fourth-order reaction diffusion equation},
R. Soc. Lond. Proc. Ser. A Math. Phys. Eng. Sci. {\bf 458} (2002), no. 2022, 1431-1446.
%
\bibitem{CEES} P Collet, J-P Eckmann, H Epstein, J Stubbe,
\textit{A global attracting set for the Kuramoto-Sivashinsky equation},
Comm. Math. Phys. {\bf 152} (1993), 203-214.
%
\bibitem{CFNT} P Constantin, C Foias, B Nicolaenko, R Temam,
\textit{Integral manifolds and inertial manifolds for dissipative partial differential equations}, 
Applied Mathematical Sciences {\bf 70},
Springer, New York, 1989.  
%
\bibitem{CFT} P Constantin, C Foias, R Temam,
{\it On the dimension of the attractors in two-dimensional turbulence},
Phys. D {\bf 30} (1988), no. 3, 284-296.
%
\bibitem{CH} M C Cross, P C Hohenberg,
{\it Pattern formation outside of  equilibrium},
Rev. Mod. Phys. {\bf 65} (1993), 851-1112.  
%
\bibitem{DFJ} R Dascaliuc, C Foias, M S Jolly,
{\it Relations between the Energy and enstrophy on the global attractor of the 2-D Navier-Stokes equations},
J. Dynam. Differential Equations {\bf 17} (2005), no. 4, 643-736.
%
\bibitem{DG} C R Doering, J D Gibbon,
{\it Applied analysis of the Navier-Stokes equations},
Cambridge Texts in Applied Mathematics,
Cambridge University Press, Cambridge, 1995.
%  
\bibitem{Go} P G\'{o}rka,
{\it Br$\acute{e}$zis-Wainger inequality on Riemannian manifolds},
J. Inequal. Appl. (2008), Art. ID 715961, 6 pages.
%
\bibitem{HN} J. M. Hyman, B. Nicolaenko,
\textit{The Kuramoto-Sivashinsky equation: a bridge between PDEs and dynamical systems},
Phys. D {\bf 18} (1986), 113-126.
%
\bibitem{Alex1} A A Ilyin,
{\it Best Constants in multiplicative inequalities for sup-norms}, 
J. London Math. Soc. {\bf 58} (1998), no. 1, 84-96.
%
\bibitem{Alex2} A A Ilyin,
\textit{Lieb-Thirring integral inequalities and their applications to attractors of Navier-Stokes equations},
Sb. Math. {\bf 196} (2005),  no. 1-2, 29-61.  
%
\bibitem{IT} A A Ilyin, E S Titi,
{\it Sharp estimates for the number of degrees  of freedom for the damped-driven 2D Navier-Stokes equations},
J. Nonlinear Sci. {\bf 16} (2006), no. 3, 233-253.        
%
\bibitem{JKT}
M Jolly, I Kevrekidis, E Titi,
\textit{Approximate inertial manifolds for the Kuramoto-Sivashinsky equation: analysis and computation},
Phys. D {\bf 44} (1990), 38-60.
%
\bibitem{KZ} L Kramer, W Zimmermann,
\textit{On the Eckhaus instability for spatially periodic patterns},
Phys. D {\bf 16} (1985), no. 2, 221-232.  
%
\bibitem{KT} Y Kuramoto, T Tsuzuki,
\textit{Persistent propagation of concentration waves in dissipative media far from thermal equilibrium},
Prog. Theor. Phys. {\bf 55} (1976), 356-369.
%
\bibitem{Mazja} V G Mazja,
{\it Sobolev spaces},
Springer-Verlag, New York, 1985.
%
\bibitem{MS} V Mazya, T Shaposhnikova,
{\it Brezis-Gallouet-Wainger type inequality for irregular domains},
Complex Var. Elliptic Equ. {\bf 56} (2011), no. 10-11, 991-1002.
%
\bibitem{PR}  L A Peletier, V Rottschafer, 
{\it Pattern selection of solutions  of the Swift-Hohenberg equation},
Phys. D {\bf 194} (2004), no. 1-2, 95-126.
%
\bibitem{PT} L A Peletier, W C Troy,
\textit{Chaotic spatial patterns described by the Extended Fisher-Kolmogorov equation}, 
J. Diff. Eqns {\bf 129} (996), 458-508.  
%
\bibitem{PW} L A Peletier, J F Williams,
{\it Some canonical bifurcations in the Swift-Hohenberg equation},
SIAM J. Appl. Dyn. Syst. {\bf 6} (2007), no. 1, 208-235. 
%
\bibitem{PM} Y Pomeau, P Manneville,
\textit{Stability and fluctuations of a spatially periodic convective flow},
J. Physique Lett. {\bf 40} (1979), no. 23, 609-612.
%
\bibitem{PZ} Y Pomeau, S Zaleski,
\textit{Wavelength selection in one-dimensional cellular structures},
J. Physique {\bf 42} (1981), no. 4, 515-528.
%
\bibitem{Robinson} J C Robinson,
{\it Infinite-dimensional dynamical systems},
Cambridge Texts in Applied Mathematics, Cambridge University Press, Cambridge, 2001.
%
\bibitem{S} G. Sivashinsky,
\textit{Nonlinear analysis of hydrodynamic instability in laminar flame. I. Derivation of basic equations},
Acta Astronautica {\bf 4} (1977), 1117-1206.
%
\bibitem{Brian} B Straughan, 
{\it The energy method, stability, and nonlinear convection}, Second edition,
Applied Mathematical Sciences {\bf 91}, Springer, New York, 2004.
%  
\bibitem{SH} J Swift, P C Hohenberg,
\textit{Hydrodynamic fluctuations at the convective instability},
Phys. Rev. A {\bf 15} (1977), 319-328.
%
\bibitem{Tadmor} E Tadmor,
\textit{The well-posedness of the Kuramoto-Sivashinsky equation},
SIAM J. Math. Anal. {\bf 17} (1986), pp. 884-893.
%
\bibitem{Temam} R Temam,
{\it Infinite-dimensional dynamical systems in mechanics and physics}, Second edition,
Applied Mathematical Sciences {\bf 68}, Springer, 1997.
%
\bibitem{X1} W Xie, 
{\it Integral representations and $L^{\infty}$ bounds for
solutions of the Helmholtz equation on arbitrary open sets in $R^2$ and $R^3$},
Differential Integral Equations {\bf 8} (1995), no. 3, 689-698.
%
\bibitem{X2} W Xie,
{\it Sharp Sobolev interpolation inequalities for the Stokes operator},
Differential Integral Equations {\bf 10} (1997), no. 2, 393-399.

\end{thebibliography}
\end{document}